\documentclass[12pt]{amsart}

\newtheorem{theorem}{Theorem}
\newtheorem{lemma}{Lemma}

\newtheorem{corollary}{Corollary}
\newtheorem{definition}{Definition}
\theoremstyle{remark}
\newtheorem{example}{Example}

\newcommand{\Cal}{\mathcal}
\newcommand{\Bb}{\mathbb}
\newcommand{\ltwod}{L^2({\mathbb R}^d)}

\begin{document}

\title{The Geometry of Sampling on Unions of Lattices}
\author{Eric Weber}
\thanks{This research was supported in part by NSF grant DMS-0200756.
}
\address{Department of Mathematics,  University of Wyoming,  Laramie, WY 82071-3036} \email{esw@uwyo.edu}
\date{\today}

\begin{abstract}
In this short note we show two results concerning sampling translation invariant subspaces of $\ltwod$ on unions of lattices.  The first result shows that the sampling transform on a union of lattices is a constant times an isometry if and only if the sampling transform on each individual lattice is so.  The second result demonstrates that the sampling transforms of two unions of lattices on two bands have orthogonal ranges if and only if correspondingly the sampling transforms of each pair of lattices have orthogonal ranges.  We then consider sampling on shifted lattices.
\end{abstract}
\maketitle




\section{Introduction}

In this note we consider sampling translation invariant subspaces on unions of lattices in $\ltwod$.  The translation invariant subspaces have the form
\[ V_E := \{ f \in \ltwod: supp( \hat{f} ) \subset E \} \]
where $E$ is a band, i.e.~a measurable subset of finite measure.  We normalize the Fourier transform so that for $f \in L^1(\Bb{R}^d) \cap \ltwod$,
\[ \hat{f}(\xi) = \int_{\Bb{R}^d} f(x) e^{-2 \pi i x \cdot \xi} dx.\]
Define, for any $y \in \Bb{R}^d$, the translation operator $T_y f(x) = f(x -y)$.  Note that $T_y V_E = V_E$ for all $y$ and all $E$.  Since $V_E$ is a reproducing kernel Hilbert space, we have that, via the inverse Fourier transform,
\begin{equation} \label{E:sample}
f(y) = \langle f, T_y \psi_E \rangle, 
\end{equation}
where $\hat{\psi}_E = \chi_E$, the indicator function of the set $E$.  Let $\{A_1, \dots, A_n \}$ be $d \times d$ invertible (real) matrices.  Define the sampling transform of the samples $A := \{A_jz + \beta_j: j=1,\dots,n; \ z \in \Bb{Z}^d\}$ by \[ \Theta_A:V_E \to \oplus_{j=1}^{n}l^2(\Bb{Z}^d): f \mapsto (f(A_1 z + \beta_1), \dots, f(A_n z + \beta_n)) \]
provided $\Theta_A$ is bounded.  A straight forward computation (see \cite{ALTW02a}) shows that $\Theta_{A}$ is bounded if and only if for each $j=1,\dots,n$, \[ \sum_{k \in \Bb{Z}^d} \chi_E(A_j^{*-1}(\xi + k)) \in L^{\infty}(\Bb{T}^d). \]
We say the samples $\{A_j z + \beta_j: j=1,\dots,n; \ z \in \Bb{Z}^d\}$ is a set of sampling for the band $E$ if $\Theta_{A}:V_E \to \oplus_{j=1}^{n}l^2(\Bb{Z}^d)$ as defined above is bounded above and below.  The samples $\{A_j z + \beta_j: j=1,\dots,n; \ z \in \Bb{Z}^d\}$ on the band $E$ are tight if $\|\Theta_{A} f\|^2 = K\|f\|^2$ for all $f \in V_E$.  The samples are an exact set of sampling for the band $E$ if the range of $\Theta_{A}$ is all of $\oplus_{j=1}^{n}l^2(\Bb{Z}^d)$.

If we oversample the space $V_E$, then the range of the sampling transform $\Theta_{A}$ will have a non-trivial orthogonal complement in $l^2(\Bb{Z}^d) \oplus \dots \oplus l^2(\Bb{Z}^d)$.  In several applications, one wishes to know the range of the sampling transform.  For example, this information is useful for denoising applications (see \cite{BT01a}).  Additionally, for multiple access communications, one also wishes to know the range of the sampling transform (see \cite{ALTW02a}).  We address this in theorems 2,3 and 4, by characterizing when two sampling transforms have orthogonal ranges.

\begin{definition}
Suppose $E$ and $F$ are bands in $\Bb{R}^d$.  Let $A_1, \dots A_n$ and $B_1,\dots, B_n$ be invertible $d \times d$ matrices, and $\{\beta_1,\dots,\beta_n,\gamma_1,\dots,\gamma_n\} \subset \Bb{R}^d$.  Let $A := \{A_j z + \beta_j: z \in \Bb{Z}^d, \ j=1,\dots,n\}$ and $B := \{B_j z + \gamma_j: z \in \Bb{Z}^d, \ j=1,\dots,n\}$ be such that the sampling transforms $\Theta_{A}$ and $\Theta_{B}$ are bounded on $E$ and $F$, respectively.  We say the samples $A$ and $B$ are orthogonal on $E$ and $F$ if $\Theta_{A} V_E$ is orthogonal to $\Theta_{B} V_F$ in $\oplus_{j=1}^{n}l^2(\Bb{Z}^d)$; equivalently $\Theta_{B}^{*} \Theta_{A} = 0$.
\end{definition}

We will use the following notation. If $C$ is an invertible matrix, denote $C^{\prime} = C^{*-1}$, where $C^*$ is the (conjugate) transpose.  If $f,g \in \ltwod$, define the $C$ bracket product to be
\[ [f,g]_C(\xi) = \sum_{m \in \Bb{Z}^d} f(\xi + Cm) \overline{g(\xi + Cm)}. \]
Let $\Pi^d$ denote the unit cube $[0,1]^d \subset \Bb{R}^d$.  Finally, define 
\[\Cal{D} = \{f \in \ltwod: \hat{f} \in L^{\infty}(\Bb{R}^d) \text{ and compactly supported}\};\]
clearly $\Cal{D}$ is a dense subspace of $\ltwod$.

\section{Unshifted Lattices}

We set out to prove the following two theorems, the main results of the paper, which relate the behavior of sampling on a union of (unshifted) lattices to the behavior of sampling on individual lattices.

\begin{theorem} \label{T:one}
Let $E$ be a band in $\Bb{R}^d$, $A_1, \dots, A_n$ be $d \times d$ matrices, and suppose that the sampling transform $\Theta_{A}$ is bounded on $V_E$.  The samples $\{A_j z: j=1,\dots,n; \ \Bb{Z}^d\}$ are tight on $E$ with constant $K$ if and only if for each $j$ the samples $\{A_j z: z \in \Bb{Z}^d\}$ are tight on $E$, with constant $K_j$.  In this case, we have $K_j = |\det A_j|^{-1}$ and $K = \sum_{j=1}^{n} K_j$.
\end{theorem}

\begin{theorem} \label{T:two}
Let $E$ and $F$ be bands in $\Bb{R}^d$, $A_1, \dots, A_n$ and $B_1,\dots,B_n$ be $d \times d$ invertible matrices, and suppose that $\Theta_{A}$ and $\Theta_B$ are bounded on $V_E$ and $V_F$, respectively.  The samples $\{A_j z: j=1,\dots,n; z \in \Bb{Z}^d\}$ and $\{B_j z:j=1,\dots,n; z \in \Bb{Z}^d\}$ are orthogonal on the bands $E$ and $F$ if and only if for $j=1,\dots,r$
\begin{equation} \label{E:8}
 \sum_{p \in \Bb{Z}^d} \chi_E(A_j^{\prime}(\xi + p)) \cdot \sum_{q \in \Bb{Z}^d} \chi_F(B_j^{\prime}(\xi + q)) = 0 \ a.e.\ \xi.
\end{equation}
Equivalently, the samples above are orthogonal if and only if the samples $\{A_jz: z \in \Bb{Z}^d\}$ and $\{B_jz: z \in \Bb{Z}^d\}$ are orthogonal on the bands $E$ and $F$ for each $j=1,\dots,r$.
\end{theorem}

As in \cite{HLW02a}, let $\Cal{P}$ be a countable index set, let $C_p$ be a $d \times d$ invertible matrix for each $p \in \Cal{P}$, and define the following:
\[ \Lambda = \cup_{p \in \Cal{P}} C_p^{\prime} \Bb{Z}^d \] 
and for $\alpha \in \Lambda$,
\[ \Cal{P}_{\alpha} = \{p \in \Cal{P}: C_p^{*} \alpha \in \Bb{Z}^d\}. \]
Note that if $\alpha = C_{p_0}^{\prime} z$ for some $z \in \Bb{Z}^d \setminus \{0\}$, then $p_0 \in \Cal{P}_{\alpha}$; if $\alpha = 0$, then $\Cal{P}_{\alpha} = \Cal{P}$.  Let $\{g_p:p \in \Cal{P}\} \subset \ltwod$.  The collection $\{T_{C_p k} g_p: p \in \Cal{P}, \ k \in \Bb{Z}^d\}$ is said to be Bessel if there exists a constant $M < \infty$ such that for all $f \in \ltwod$,
\[ \sum_{p \in \Cal{P}} \sum_{k \in \Bb{Z}^d} | \langle f, T_{C_p k} g_p \rangle |^2 \leq M \|f \|^2. \]
The collection $\{T_{C_p k} g_p: p \in \Cal{P}, \ k \in \Bb{Z}^d\}$ satisfies the local boundedness condition if for every $f \in \Cal{D}$,
\[ L(f) := \sum_{p \in \Cal{P}} \sum_{k \in \Bb{Z}^d} \int_{supp \hat{f}} |\hat{f} (\xi + C_{p}^{\prime} k) |^2 |\det C_p|^{-1} |\hat{g}_{p}(\xi)|^2 d \xi < \infty. \]

\begin{lemma} \label{L:one}
Suppose the sets $G:=\{T_{C_p k} g_p: p \in \Cal{P}, \ k \in \Bb{Z}^d\}$ and $H := \{T_{C_p k} h_p: p \in \Cal{P}, \ k \in \Bb{Z}^d\}$ are both Bessel and both satisfy the local boundedness condition $L(f) < \infty$ for all $f \in \Cal{D}$, as defined above.  Define the operator
\[ \Omega_{H,G}: \ltwod \to \ltwod: f \mapsto \sum_{p \in \Cal{P}} \sum_{k \in \Bb{Z}^d} \langle f, T_{C_p k} g_p \rangle T_{C_p k} h_p. \]
The operator $\Omega_{H,G}$ commutes with the translation operators $T_y$ for all $y \in \Bb{R}^d$ if and only if for all $\alpha \in \Lambda \setminus \{0\}$,
\[ \sum_{p \in \Cal{P}_{\alpha}} |\det C_p|^{-1} \overline{\hat{g}_{p}(\xi)} \hat{h}_{p}(\xi + \alpha) = 0 \ a.e. \ \xi.\]
In this case, $\Omega_{H,G}$ is a Fourier multiplier with symbol
\[ s(\xi) := \sum_{p \in \Cal{P}} |\det C_p|^{-1} \overline{\hat{g}_{p}(\xi)}\hat{h}_{p}(\xi), \]
i.e.~$\widehat{\Omega_{H,G} f}(\xi) = s(\xi) \hat{f}(\xi)$.
\end{lemma}
\begin{proof}
For $f \in \Cal{D}$, define the continuous function
\[ w_f(x) = \langle \Omega_{H,G} T_x f , T_x f \rangle. \]
If $\Omega_{H,G}$ commutes with all $T_x$ for $x \in \Bb{R}^d$, then clearly $w_f(x)$ is constant for all $f \in \Cal{D}$.  Conversely, if $w_f(x)$ is constant for all $f \in \Cal{D}$, then $\langle T_{-x} \Omega_{H,G} T_x f, f \rangle = \langle \Omega_{H,G} f, f \rangle$, whence by the polarization identity, $T_{-x} \Omega_{H,G} T_x = \Omega_{H,G}$, and thus $\Omega_{H,G} T_x = T_x \Omega_{H,G}$.

By \cite[Proposition 2.4]{HLW02a}, $w_f(x)$ coincides pointwise with the almost periodic function 
\[\sum_{\alpha \in \Lambda} \hat{w}_f(\alpha) e^{2 \pi i \alpha \cdot x},\] 
where
\[ \hat{w}_f(\alpha) = \int_{\Bb{R}^d} \hat{f}(\xi) \overline{\hat{f}(\xi + \alpha)} \sum_{p \in \Cal{P}_{\alpha}} |\det C_p|^{-1} \overline{\hat{g}_p(\xi)}\hat{h}_p(\xi + \alpha) d \xi.\]
By \cite[Lemma 2.5]{HLW02a} and the proof of Theorem 2.1 in \cite{HLW02a}, $w_f(x)$ is constant for all $f \in \Cal{D}$ if and only if for all $\alpha \in \Lambda \setminus \{0\}$,
\[ \sum_{p \in \Cal{P}_{\alpha}} |\det C_p|^{-1} \overline{\hat{g}_p(\xi)}\hat{h}_p(\xi + \alpha) = 0 a.e. \ \xi.\]

It is well known that if $\Omega_{H,G}$ commutes with $T_y$ for all $y \in \Bb{R}^d$, then it is a Fourier multiplier.  Evaluating $w_f(x)$ at $x = 0$ yields
\[ w_f(0) = \int_{\Bb{R}^d} \hat{f}(\xi) \overline{\hat{f}(\xi)} \sum_{p \in \Cal{P}} |\det C_p|^{-1} \overline{\hat{g}_{p}(\xi)}\hat{h}_{p}(\xi) d \xi = \langle \Omega_{H,G} f, f \rangle.\]
Therefore, the symbol of $\Omega_{H,G}$ is $s(\xi)$ as above.
\end{proof}

\begin{lemma} \label{L:two}
Suppose $\{g_p\}_{p \in \Cal{P}} \subset V_E$ and let $G$ be as in lemma \ref{L:one}.  Define $\Theta_{G}$
\[ \Theta_{G}: \ltwod \to l^2(\Bb{Z}^d \times \Cal{P}): f \mapsto (\langle f, T_{C_p k} g_p \rangle)_{(k, p)}. \]
We have $\|\Theta_{G} g \|^2  = K \|g \|^2$ for all $g \in V_{E}$ if and only if for all $\alpha \in \Lambda$,
\begin{equation} \label{E:P_q}
\sum_{p \in \Cal{P}_{\alpha}} |\det C_p|^{-1} \overline{\hat{g}_{p}(\xi)} \hat{g}_{p}(\xi + \alpha) = K \delta_{\alpha, 0} \chi_{E}(\xi) \ a.e. \ \xi,
\end{equation}
where $\delta_{\alpha, 0}$ is the Kronecker delta.
\end{lemma}

\begin{proof}
Note that, as defined in lemma \ref{L:one}, $\Omega_{G,G} = \Theta^{*}_{G} \Theta_{G}$, and hence $\|\Theta_{G} g\|^2 = \langle \Omega_{G,G}g, g \rangle$.

Suppose equation \ref{E:P_q} holds.  Then, by lemma \ref{L:one}, $\Omega_{G,G}$ is a Fourier multiplier, with symbol
\[ s(\xi) = \sum_{p \in \Cal{P}} |\det C_p|^{-1} |\hat{g}_p(\xi)|^{2} = K \chi_{E}(\xi). \]
Thus, 
\[ \|\Theta_{G} g\|^2 = \langle \Omega_{G,G} g , g \rangle = \langle s(\xi) \hat{g}(\xi), \hat{g}(\xi) \rangle = \langle K g, g \rangle = K \|g\|^2. \]

Conversely, suppose $\|\Theta_{G} g \|^2 = K \|g\|^2$ for all $g \in V_{E}$.  Then, by the above computation, it follows that $\langle \Omega_{G,G} g, g \rangle = \langle K g, g \rangle$, whence by the polarization identity, $\Omega_{G,G}$ is a Fourier multiplier with symbol $s(\xi) = K \ a.e. \ \xi$.  It now follows that equation \ref{E:P_q} holds.
\end{proof}

We wish to apply lemmas \ref{L:one} and \ref{L:two} to theorems \ref{T:one} and \ref{T:two}.  Let $A:=\{A_j z + \beta_j:j=1,\dots,n; \ z \in \Bb{Z}^d\}$, where $A_j$ is an invertible matrix for each $j$, and let $E$ be a band such that $\Theta_{A}$ is bounded on $V_E$.  By equation \ref{E:sample}, the sampling transform $\Theta_{A}$ above can be written in terms of the collection
\begin{equation}  \label{E:coll}
\{ T_{A_j z} T_{\beta_j} \psi_{E}: j=1,\dots,n; \ z \in \Bb{Z}^d \},
\end{equation}
as follows:
\[ ( f(A_1 z + \beta_1), \dots, f(A_n z + \beta_n)) = (\langle f, T_{A_1 z} T_{\beta_1} \psi_{E} \rangle, \dots, \langle f, T_{A_n z} T_{\beta_n} \psi_{E} \rangle ). \]
Hence, questions of tightness and orthogonality of sampling can be stated in terms of the collection (\ref{E:coll}).  In order to apply the lemmas, we need to verify that (\ref{E:coll}) satisfies the local boundedness condition.

\begin{lemma} \label{L:three}
Let $A$, $E$ be as above.  If the sampling transform $\Theta_A$ is bounded on $V_E$, then for all $f \in \Cal{D}$, we have
\[ L(f) = \sum_{j=1}^{n} \sum_{m \in \Bb{Z}^d} \int_{supp \hat{f}} |\hat{f}(\xi + A_j^{\prime}m)|^2 |\det A_j|^{-1} |\chi_{E}(\xi)|^2 d \xi < \infty.\]
\end{lemma}
\begin{proof}
Let $Q = \|f\|_{\infty}$.  Since $\Theta_{A}$ is bounded on $V_E$, then the cardinality of the set
\[ \{m \in \Bb{Z}^d: \xi + A_j^{\prime}m \in E \} \]
is essentially bounded in $\xi$.  Let this bound be $R$.  Clearly, we have
\[ L(f) \leq \sum_{j=1}^{n} |\det A_j|^{-1}\int_{E} Q R \ d \xi < \infty. \] 
\end{proof}

The following lemma is a variation of a lemma in \cite{HLW02a}; it will be needed for theorem \ref{T:two}.
\begin{lemma} \label{L:four}
Let $G$ be as in lemma \ref{L:one}; let $D_p$ be a $d \times d$ invertible matrix for each $p \in \Cal{P}$, and suppose $H:=\{T_{D_p k} h_p: p \in \Cal{P}, \ k \in \Bb{Z}^d\}$ is Bessel and satisfies the local boundedness condition.  Define $\Omega_{H,G}$ analogous to lemma \ref{L:one}:
\[ \Omega_{H,G}: \ltwod \to \ltwod: f \mapsto \sum_{p \in \Cal{P}} \sum_{k \in \Bb{Z}^d} \langle f, T_{C_p k} g_p \rangle T_{D_p k} h_p. \]
For all $f,g \in \Cal{D}$,
\[ \langle \Omega_{H,G} f, g \rangle = \sum_{p \in \Cal{P}} \frac{1}{|\det C_p D_p |} \int_{\Pi^d} [\hat{f}, \hat{g}_p]_{C_p^{\prime}}(C_p^{\prime} \xi) [\hat{h}_p, \hat{g}]_{D_p^{\prime}} (D_p^{\prime} \xi) d \xi. \]
\end{lemma}
\begin{proof}
We have via the Fourier transform
\begin{align*} 
\langle \Omega_{H,G} f, g \rangle &= \sum_{p \in \Cal{P}} \sum_{k \in \Bb{Z}^d} \langle f, T_{C_pk}g_p \rangle \langle T_{D_pk}h_p, g \rangle\\
&= \sum_{p \in \Cal{P}} \sum_{k \in \Bb{Z}^d} \int_{\Bb{R}^d} \hat{f}(\xi) e^{2\pi i C_p k \cdot \xi} \overline{\hat{g}_p(\xi)} d \xi \cdot \int_{\Bb{R}^d} e^{-2\pi i D_p k \cdot \omega} \hat{h}_p(\omega) \overline{\hat{g}(\omega)} d \omega.
\end{align*}
Each of the above integrals can be written as follows, by variable substitutions.
\begin{align*}
\int_{\Bb{R}^d} \hat{f}(\xi) e^{2\pi i C_p k \cdot \xi} \overline{\hat{g}_p(\xi)} d \xi &= |\det C_p|^{-1} \int_{\Bb{R}^d} \hat{f}(C_p^{\prime}\xi) e^{2\pi i k \cdot \xi} \overline{\hat{g}_p(C_p^{\prime}\xi)} d \xi \\
&= \sum_{m \in \Bb{Z}^d} |\det C_p|^{-1} \int_{\Pi^d} \hat{f}(C_p^{\prime}(\xi + m)) e^{2\pi i k \cdot \xi} \overline{\hat{g}_p(C_p^{\prime}(\xi + m))} d \xi \\
&= |\det C_p|^{-1} \int_{\Pi^d} [\hat{f}, \hat{g}_p]_{C_p^{\prime}}(C_p^{\prime} \xi) e^{2 \pi i k \cdot \xi} d \xi.
\end{align*}
The interchanging of the sum and integral is valid since $f,g \in \Cal{D}$.  Therefore, we have
\begin{equation} \label{E:sum}
\langle \Omega_{H,G} f , g \rangle = \sum_{p \in \Cal{P}}  \frac{1}{|\det C_p D_p|}  \sum_{k \in \Bb{Z}^d} \int_{\Pi^d} [\hat{f}, \hat{g}_p]_{C_p^{\prime}}(C_p^{\prime} \xi) e^{2 \pi i k \cdot \xi} d \xi \cdot \int_{\Pi^d} [\hat{h}_p, \hat{g}]_{D_p^{\prime}}(D_p^{\prime} \omega) e^{-2 \pi i k \cdot \omega} d \omega 
\end{equation}
where the sum over $k$ is really the inner product of the Fourier coefficients of the two functions
\[ \overline{[\hat{f}, \hat{g}_p]_{C_p^{\prime}}(C_p^{\prime} \xi)} \text{ and } [\hat{h}_p, \hat{g}]_{D_p^{\prime}}(D_p^{\prime} \omega). \]
Thus, the sum (\ref{E:sum}) is
\[ \sum_{p \in \Cal{P}} \frac{1}{|\det C_p D_p |} \int_{\Pi^d} [\hat{f}, \hat{g}_p]_{C_p^{\prime}}(C_p^{\prime} \xi) [\hat{h}_p, \hat{g}]_{D_p^{\prime}}(D_p^{\prime} \xi) d \xi. \]
\end{proof}

\begin{proof}[Proof of Theorem 1.]
Let $\Cal{P} = \{1,\dots,n\}$, $C_j = A_j$, and $g_j = \psi_E$ for $j=1,\dots,n$.  Thus, $\Lambda = \cup_{j=1}^{n} A_j^{\prime} \Bb{Z}^d$ and for $\alpha \in \Lambda$, $\Cal{P}_{\alpha} = \{ j \in \{1,\dots,n\}: A_j^{*} \alpha \in \Bb{Z}^d\}$ as above.  Note that if $\alpha = A_{j_0}^{\prime} z$ for some $z \in \Bb{Z}^d \setminus \{0\}$, then $j_0 \in \Cal{P}_{\alpha}$.  By equation \ref{E:sample} and lemmas \ref{L:two} and \ref{L:three}, we have that
\[ \|\Theta_{A}f\|^2 = K \|f\|^2 \]
if and only if
\[ \sum_{j \in \Cal{P}_\alpha} |\det A_j|^{-1} \chi_E(\xi)\chi_E(\xi + \alpha) = K \delta_{\alpha,0} \chi_E(\xi) \]
for all $\alpha \in \Lambda$. 

It follows that for $z \in \Bb{Z}^d \setminus \{0\}$ and for each $j \in \{1,\dots,n\}$, 
\[ \chi_E(\xi) \chi_E(\xi + A_j^{\prime}z) = 0 \ a.e.\ \xi.\]
Therefore, for each $j$,
\[ \sum_{m \in \Bb{Z}^d} \chi_E(A_j^{\prime}(\xi + m)) \leq 1 \ a.e.\ \xi.\]
Since the sum is bounded above by 1, for all $f \in V_E \cap \Cal{D}$, 
\[ \|\Theta_{A_j}f\|^2 = |\det A_j|^{-1} \|f\|^2, \]
(see \cite{K65a} and \cite{ALTW02a}).
\end{proof}

\begin{proof}[Proof of Theorem 2.]
In \cite{ALTW02a}, we show that equation \ref{E:8} is sufficient for $\Theta_{A}V_E \perp \Theta_{B}V_F$; thus we wish to show that it is necessary.  Let $E_0 \subset E$ and $F_0 \subset F$ be any measurable bounded sets.  Define $g_0 \in V_E$ and $f_0 \in V_F$ by $\hat{g}_0 = \chi_{E_0}$ and $\hat{f}_0 = \chi_{F_0}$.  Suppose $\Theta_{A}V_E \perp \Theta_{B}V_F$; by lemma \ref{L:four} we have
\[\sum_{j=1}^{n} |\det A_j B_j|^{-1} \int_{\Pi^d} [\chi_{E_0}, \chi_E]_{A_j^{\prime}}(A_j^{\prime}\xi) [\chi_F, \chi_{F_0}]_{B_j^{\prime}}(B_j^{\prime} \xi) = 0. \]
Therefore, for each $j=1, \dots, r$,
\[ \int_{\Pi^d} [\chi_{E_0}, \chi_E]_{A_j^{\prime}}(A_j^{\prime}\xi) [\chi_F, \chi_{F_0}]_{B_j^{\prime}}(B_j^{\prime} \xi) = 0. \]
However, expanding the above integrand yields
\[ \int_{\Pi^d} \sum_{p \in \Bb{Z}^d} \chi_{E_0}(A_j^{\prime}(\xi + p)) \sum_{q \in \Bb{Z}^d} \chi_{F_0}(B_j^{\prime}(\xi + q)) = 0.\]
The result now follows since $E_0$ and $F_0$ were arbitrary measurable bounded subsets.
\end{proof}

\section{Shifted Lattices}

We begin with the following remark.  Both theorems 1 and 2 are no longer valid if the lattices are allowed to be shifted.  As a trivial example for theorem 1, consider the canonical bandlimited function space $V_{[-1/2, 1/2]}$ and the lattices $2\Bb{Z}$ and $2\Bb{Z} + 1$; individually neither form sets of sampling, but together they form an exact set of sampling.  An example corresponding to theorem 2 is as follows.

\begin{example}
Let $E=[0,1]$, $F=[-1,0]$ and define the following sampling transforms on $V_E$ and $V_F$ respectively:
\begin{align*}
\Theta_{A}&: V_E \to l^2(\Bb{Z}) \oplus l^2(\Bb{Z}): f \mapsto (f(n), f(n + 1/2)) \\
\Theta_{B}&: V_F \to l^2(\Bb{Z}) \oplus l^2(\Bb{Z}): g \mapsto (g(n + 1/2), g(n)).
\end{align*}
Define $\Theta_{A_j}$ and $\Theta_{B_j}$, for $j=1,2$, to be the first and second coordinates of $\Theta_{A}$ and $\Theta_{B}$ respectively.  It is clear that $\Theta_{A_j} V_E = \Theta_{B_j} V_F = l^2(\Bb{Z})$.  However, we claim that the ranges of $\Theta_{A}$ and $\Theta_{B}$ are orthogonal.

Define a unitary operator $U:V_E \to V_F: \hat{f}(\xi) \mapsto \hat{f}(\xi + 1)$.  By the polarization identity, it suffices to show that
\[ \langle \Theta_{A} f, \Theta_{B} U f \rangle = 0 \]
for all $f \in V_E$.  Note that 
\[ \langle \Theta_{A} f, \Theta_{B} U f \rangle = \langle \Theta_{A_1} f, \Theta_{B_1} U f \rangle + \langle \Theta_{A_2} f, \Theta_{B_2} U f \rangle.\]
Note also that 
\begin{align*}
[Uf](x) &= \int_{-1}^{0} \widehat{Uf}(\xi) e^{-2\pi i x \xi} d \xi \\
&= \int_{0}^{1} e^{2 \pi i x} \hat{f}(\xi) e^{-2\pi i x \xi} d \xi
\end{align*}
via the variable substitution $\xi \to \xi - 1$.

Therefore, we have that 
\begin{align*}
\langle \Theta_{A_1} f, \Theta_{B_1} U f \rangle &= \sum_{n \in \Bb{Z}} \int_{0}^{1} \hat{f}(\xi) e^{-2 \pi i n \xi} d \xi \cdot \overline{\int_{0}^{1} e^{2 \pi i (n + 1/2)} \hat{f}(\xi) e^{-2 \pi i (n + 1/2) \xi} d \xi} \\
&= \sum_{n \in \Bb{Z}} \int_{0}^{1} \hat{f}(\xi) e^{-2 \pi i n \xi} d \xi \cdot \overline{\int_{0}^{1} - e^{\pi i \xi} \hat{f}(\xi) e^{-2 \pi i n \xi} d \xi} \\&= \langle \hat{f}(\xi), - e^{\pi i \xi} \hat{f}(\xi) \rangle_{L^2{[0,1]}}.
\end{align*}
Likewise,
\begin{align*}
\langle \Theta_{A_2} f, \Theta_{B_2} U f \rangle &= \sum_{n \in \Bb{Z}} \int_{0}^{1} \hat{f}(\xi) e^{-2 \pi i (n + 1/2) \xi} d \xi \cdot \overline{\int_{0}^{1} e^{2 \pi i n} \hat{f}(\xi) e^{-2 \pi i n \xi} d \xi} \\
&= \sum_{n \in \Bb{Z}} \int_{0}^{1} e^{- \pi i \xi} \hat{f}(\xi) e^{-2 \pi i n \xi} d \xi \cdot \overline{\int_{0}^{1} \hat{f}(\xi) e^{-2 \pi i n \xi} d \xi} \\
&= \langle e^{- \pi i \xi} \hat{f}(\xi), \hat{f}(\xi) \rangle_{L^2{[0,1]}}.
\end{align*}
Therefore,
\[ \langle \Theta_{A} f, \Theta_{B} U f \rangle = \langle \hat{f}(\xi), - e^{\pi i \xi} \hat{f}(\xi) \rangle + \langle e^{- \pi i \xi} \hat{f}(\xi), \hat{f}(\xi) \rangle = 0.\]
\end{example}

Let $m$ denote Lebesgue measure on $\Bb{R}^d$.
\begin{corollary}
The samples $\{A_j z + \beta_j:j=1,\dots,n; \ z \in \Bb{Z}^d\}$ are tight on the band $E$ if and only if for each $\alpha \in \Lambda \setminus \{0\}$, either
\begin{enumerate}
\item $m(E \cap (E + \alpha)) = 0$, or
\item $\sum_{j \in \Cal{P}_{\alpha}} |det A_j|^{-1} e^{-2 \pi i \beta_j \cdot \alpha} = 0.$
\end{enumerate}
\end{corollary}

\begin{proof}
Let $\Cal{P}$ and $C_j$ be as in the proof of Theorem \ref{T:one}; let $g_j = T_{\beta_j} \psi_{E}$.  The samples $\{A_j z + \beta_j:j=1,\dots,n; \ z \in \Bb{Z}^d\}$ are tight on the band $E$ if and only if equation \ref{E:P_q} in Lemma \ref{L:two} is satisfied.  Note that the condition for $\alpha = 0$ is automatically satisfied, with $K = \sum_{j=1}^{n} |\det A_j|^{-1}$.  For $\alpha \neq 0$, we must have
\begin{align*}
\sum_{j \in \Cal{P}_{\alpha}} \overline{\hat{g}_j(\xi)}\hat{g_j}(\xi + \alpha) &=  \sum_{j \in \Cal{P}_{\alpha}} |\det A_j|^{-1} e^{2 \pi i \beta_j \cdot \xi} \chi_{E}(\xi) e^{-2 \pi i \beta_j \cdot (\xi + \alpha)}\chi_{E}(\xi + \alpha) \\
&= \sum_{j \in \Cal{P}_{\alpha}} |\det A_j|^{-1} e^{-2 \pi i \beta_j \cdot \alpha} \cdot  \chi_{E}(\xi) \chi_{E}(\xi + \alpha) \\
&= 0,
\end{align*}
from which the statement follows.
\end{proof}

An immediate corollary of lemma 2 regarding orthogonality is as follows.
\begin{corollary} \label{C:two}
The samples $\{A_j z + \beta_j:j=1,\dots,n; \ z \in \Bb{Z}^d\}$ and $\{B_j z + \gamma_j:j=1,\dots,n; \ z \in \Bb{Z}^d\}$ are orthogonal on the bands $E$ and $F$ if and only if for all $f \in V_{E \cup F} \cap \Cal{D}$
\[ \sum_{j=1}^{n} |\det A_j B_j|^{-1} \int_{\Pi^d} [\hat{f}, e^{-2 \pi i \beta_j(\cdot)} \chi_{E}]_{A_j^{\prime}}(A_j^{\prime}\xi) [e^{-2 \pi i \gamma_j(\cdot)} \chi_{F}, \hat{f}]_{B_j^{\prime}}(B_j^{\prime} \xi) d \xi = 0.\]
\end{corollary}
A (better) complete characterization of determining orthogonality of shifted lattices on bands $E$ and $F$ is out of reach at this time.  Note that in lemma \ref{L:one}, $C_p$ appears in both sequences.  We wish to consider orthogonal samples with different matrices; the techniques used to prove lemma \ref{L:one} do not readily extend to this case.  We conclude, however, with the following two special cases: 1) $B_j = A_j$, and 2) $A_j = A$ and $B_j = B$.  

\begin{lemma} \label{L:five}
Suppose $G$ and $H$ and $\Omega_{H,G}$ are as in lemma \ref{L:one}; $\Omega_{H,G} = 0$ if and only if for all $\alpha \in \Lambda$,
\begin{equation} \label{E:zero}
 \sum_{p \in \Cal{P}_{\alpha}} |\det C_p|^{-1} \overline{\hat{g}_{p}(\xi)} \hat{h}_{p}(\xi + \alpha) = 0 \ a.e. \ \xi.
\end{equation}
\end{lemma}
\begin{proof}
If equation \ref{E:zero} holds for all $\alpha \in \Lambda \setminus \{0\}$, then $\Omega_{H,G}$ is a Fourier multiplier; if the equation also holds for $\alpha = 0$, then the symbol of $\Omega_{H,G}$ is $0$.  Conversely, if $\Omega_{H,G} = 0$, then $\Omega_{H,G}$ commutes with all translation operators $T_y$, hence equation \ref{E:zero} holds for $\alpha \in \Lambda \setminus \{0\}$; moreover, since $\Omega_{H,G} = 0$, then the equation must also hold for $\alpha = 0$.
\end{proof}

Let $\Lambda$ be as above; for $q \in \Bb{R}^d$, define $I_{q,\alpha} := \{ j \in \Cal{P}_{\alpha}: \gamma_j - \beta_j = q\}$.

\begin{theorem} \label{T:Three}
The samples $A:=\{A_j z + \beta_j:i=1,\dots,n; \ z \in \Bb{Z}^d \}$ and $B:=\{A_j z + \gamma_j:j=1,\dots,n; z \in \Bb{Z}^d \}$ are orthogonal on the bands $E$ and $F$, respectively, if and only if for all $\alpha \in \Lambda$ either
\begin{enumerate}
\item $m(E \cap (F + \alpha)) = 0$, or
\item $\sum_{j \in I_{q, \alpha}} |\det A_j|^{-1} e^{-2 \pi i \gamma_j \cdot \alpha} = 0,$ for all $q \in \Bb{R}^d$.
\end{enumerate}
In particular, $m(E \cap F) = 0$.
\end{theorem}
\begin{proof}
As above, let $\Cal{P}$, $C_j$, and $g_j$ be as in Corollary \ref{C:two}; let $h_j = T_{\gamma_j}\psi_{F}$.  The samples $A$ and $B$ are orthogonal on $E$ and $F$ respectively if and only if $\Theta_{B}^{*} \Theta_{A} = 0$.  Let $G:=\{T_{A_j k} g_j: j=1,\dots,n; \ k \in \Bb{Z}^d\}$ and $H:=\{T_{A_j k} h_j: j=1,\dots,n; \ k \in \Bb{Z}^d\}$.  Note that $\Omega_{H,G} = 0$ if and only if $\Theta_{B}^{*} \Theta_{A} = 0$.

Applying Lemma \ref{L:five} yields that for all $\alpha \in \Lambda$, we must have
\[ \sum_{j \in \Cal{P}_{\alpha}} |\det A_j|^{-1} e^{2 \pi i \beta_j \cdot \xi} \chi_{E}(\xi) e^{-2 \pi i \gamma_j \cdot (\xi + \alpha)} \chi_{F}(\xi + \alpha) = 0. \]
Clearly, if $m(E \cap (F + \alpha)) \neq 0$, then the trigonometric polynomial
\[ \sum_{q = \gamma_j - \beta_j} \sum_{k \in I_{q,\alpha}} |\det A_k|^{-1} e^{- 2 \pi i \gamma_k \cdot \alpha} e^{-2 \pi i q \cdot \xi} = 0.\]
The converse follows directly.
\end{proof}

For $q \in \Bb{R}^d$, define $J_q := \{ j \in \{1, \dots,n\}: B^{-1} \gamma_j - A^{-1} \beta_j = q\}$.

\begin{theorem}
The samples $A:=\{ A z + \beta_j:j=1,\dots,n; z \in \Bb{Z}^d \}$ and $B:=\{ B z + \gamma_j: j=1,\dots,n; z \in \Bb{Z}^d\}$ on the bands $E$ and $F$, respectively, are orthogonal if and only if for all $m \in \Bb{Z}^d$, either
\begin{enumerate}
\item $m(A^{*}E \cap (B^{*}F + m)) = 0$, or
\item $\sum_{j \in J_q} e^{-2 \pi i \gamma_j B^{\prime} m} = 0,$ for all $q \in \Bb{R}^d$.
\end{enumerate}
In particular, $m(A^{*}E \cap B^{*}F) = 0$.
\end{theorem} 
\begin{proof}
The proof follows the proof of Theorem \ref{T:three} with the following modification.  If $C$ is a $d \times d$ invertible matrix, then define the unitary dilation operator $D_C$ on $\ltwod$ by $D_C g(x) = \sqrt{|\det C|} g(C x)$.  We have
\begin{align*} 
\sum_{j=1}^{n} \sum_{k \in \Bb{Z}^d} \langle f, T_{A k} T_{\beta_j} \psi_{E} \rangle & \langle T_{B k} T_{\gamma_j} \psi_{F}, g \rangle \\
&= \sum_{j=1}^{n} \sum_{k \in \Bb{Z}^d} \langle D_A f, D_A T_{A k} T_{\beta_j} \psi_{E} \rangle \langle D_B T_{B k} T_{\gamma_j} \psi_{F}, D_B g \rangle \\
&= \sum_{j=1}^{n} \sum_{k \in \Bb{Z}^d} \langle D_A f, T_k D_A T_{\beta_j} \psi_{E} \rangle \langle T_k D_B T_{\gamma_j} \psi_{F}, D_B g \rangle \\
&= \langle \Omega_{H,G} D_{A} f, D_{B} g \rangle
\end{align*}
where $G:=\{T_k D_{A} T_{\beta_j} \psi_{E}: j=1,\dots,n; \ k \in \Bb{Z}^d$ and $H:= \{T_k D_{B} T_{\gamma_j} \psi_{F}: j=1,\dots,n; \ k \in \Bb{Z}^d$.  It follows, since $D_A$ and $D_B$ are unitary, that the samples $A$ and $B$ are orthogonal on the bands $E$ and $F$ if and only if $\Omega_{H,G} = 0$.
Thus, we apply Lemma \ref{L:five} with $\Lambda = \Bb{Z}^d$, $g_j = D_A T_{\beta_j} \psi_E$, and $h_j = D_B T_{\gamma_j} \psi_F$.  The computations are as in Theorem \ref{T:three}.
\end{proof}

\emph{Acknowledgement.} We thank the anonymous referee for comments and suggestions which greatly improved the presentation of the paper.

\bibliographystyle{amsplain}
\bibliography{gsl}

\providecommand{\bysame}{\leavevmode\hbox to3em{\hrulefill}\thinspace}
\providecommand{\MR}{\relax\ifhmode\unskip\space\fi MR }
\providecommand{\MRhref}[2]{%
  \href{http://www.ams.org/mathscinet-getitem?mr=#1}{#2}
}
\providecommand{\href}[2]{#2}
\begin{thebibliography}{1}

\bibitem{ALTW02a}
A.~Aldroubi, D.~Larson, W.~S. Tang, and E.~Weber, \emph{The geometry of frame
  representations of abelian groups}, preprint, posted on arXiv.org, math.FA/0308250.

\bibitem{BF01b}
H.~Behmard and A.~Faridani, \emph{Sampling of bandlimited functions on unions
  of shifted lattices}, J. Fourier Anal. Appl. \textbf{8} (2001), no.~1,
  43--58.

\bibitem{BT01a}
J.~Benedetto and O.~Treiber, \emph{Wavelet frames: multiresolution analysis and
  extension principles}, Wavelet transforms and time-frequency signal analysis,
  Appl. Numer. Harmon. Anal., Birkhauser Boston, 2001, pp.~3--36.

\bibitem{HLW02a}
E.~Hernandez, D.~Labate, and G.~Weiss, \emph{A unified characterization of
  reproducing systems generated by a finite family {II}}, J. Geom. Anal.
  \textbf{12} (2002), no.~4, 615--662.

\bibitem{K65a}
I.~Kluv\'anek, \emph{Sampling theory in abstract harmonic analysis}, Mat.-Fyz.
  Capopis Sloven. Akad. Vied. \textbf{15} (1965), 43--48.

\bibitem{Wal96a}
David Walnut, \emph{Nonperiodic sampling of bandlimited functions on unions of
  rectangular lattices}, J. Fourier Anal. Appl. \textbf{2} (1996), no.~5,
  435--452.

\end{thebibliography}
\nocite{Wal96a}
\nocite{BF01b}

\end{document}